# Gauss-Type Quadrature Rules Based on Identity-Type Functions

**M. A. Bokhari**[1] **and Asghar Qadir**[1,2]

[1] Department of Mathematics and Statistics
King Fahd Univ. of Petroleum & Minerals
Dhahran, Saudi Arabia
mbokhari@kfupm.edu.sa

[2] Center for Applied Mathematics & Physics,
National University of Science & Technology
Rawalpindi, Pakistan
aqadirmath@yahoo.com

*Some Gauss-type quadrature rules over [0, 1], which involve values and/or the derivative of the integrand at 0 and/or 1, are investigated. Our work is based on the orthogonal polynomials with respect to the linear weight function $\omega(t) := 1 - t$ over [0, 1], which belong to the class of Jacobi polynomials. It has been shown that the $\omega(t)$-orthogonal polynomials arise from the recently developed "identity-type functions". Along the lines of Golub's work, the nodes and weights of the quadrature rules are computed from Jacobi-type matrices. Based on the properties of the identity-type functions, we discover that the resulting matrices consist of simple rational entries. Computational procedures for the derived rules are tested on different integrands. The proposed methods have some advantage over the respective Gauss-type rules with respect to the constant weight function $\omega(t) := 1$ over [0, 1].*

**Keywords:** Gauss / Gauss-Radau / Gauss-Lobatto quadrature rules, Jacobi-matrix, Hypergeometric series, Identity-type polynomials, 3-term recurrence relation, Recursion coefficients.

## Introduction

Numerical integration is usually utilized when analytical techniques fail. On the other hand, when a high degree of accuracy is required, the number of steps related to a quadrature rule can become too large, creating inaccuracies due to error propagation. For this reason one looks for efficient schemes that can reduce the number of steps while retaining the accuracy at the desired level. The Gaussian quadrature rule is one such scheme. It requires the weighted integrand's values at the zeros of orthogonal polynomials, which lie inside the interval of integration, say [0, 1]. Several other schemes exist that deal in addition with the values of the integrand and/or its derivatives at 0 and/or 1.

Among these, Gauss-Radau and Gauss-Lobatto rules are quite prominent. These rules also prove useful in determining approximate solutions of boundary value problems. A lot of work is still going on in the modification of quadrature rules and their computational aspects. Gautschi and Li discussed some formulae with double end-point nodes for all four Chebyshev weight functions in [5]. These are referred to as generalized Gauss-Radau and Gauss-Lobatto rules in the literature. In 1983, Golub and Kautsky [8] derived algorithms for the evaluation of Gauss knots in the presence of fixed knots by modifying a Jacobi matrix. Gautschi developed computational methods for generating Gauss-type quadrature formulae having nodes of arbitrary multiplicity at one or both end points of the interval of integration [3]. Li studied Kronrod extensions of Gauss-Radau and Gauss-Lobatto formulae having end points of multiplicity 2 and derived explicit formulae associated with the end points for the Chebyshev weights [10].

In a totally different context [1], a class of "identity-type functions" based on Gauss hypergeometric functions has been recently defined, which has interesting properties. Specifically, it includes what are called "identity-type polynomials" [1, (3.1)] that satisfy an orthogonality property with respect to weight function $\omega(t) := 1 - t$. These polynomials give rise to "Jacobi-type" polynomials which may be used to derive Gauss-type quadrature rules which involve values and/or the derivative of the integrand at 0 and/or 1. However, the recursion coefficients directly evaluated from the formulae (5) are not appropriate due to accumulation of round off errors in the calculation of nodes and weights of the Quadrature rules via Jacobi matrix (6). Nevertheless, certain properties of the identity-type functions enabled us to derive explicit representation of these coefficients in the form of simple rational sequences. Such representation which guarantees high accuracy of Gauss-type quadrature rules is not available in the literature.

Since hypergeometric functions arise in the solution of Sturm-Liouville problems, the derived polynomials could be applicable to the quadrature schemes, especially if the quadrature is used for solving differential equations. The objective of the present work is to examine the structure and computational aspects of Gauss-type quadrature rules relevant to these polynomials.

The plan of the paper is as follows. In the next section we set up the Gauss, Gauss-Radau and Gauss-Lobatto schemes subject to weight function $\omega(t) := 1 - t$ over [0, 1]. In the subsequent sections we provide a review of the identity-type polynomials and then discuss orthogonality of their specific factors. We also provide, in the same section, a simple expression for the recursion coefficients required to compute nodes and weights of the Quadrature rules. Section 5 deals with some numerical examples. A brief summary and discussion of our results will be given in the last section.

# Gauss-Type Quadrature Procedures

We begin with setting $f_\omega(t) := \begin{cases} \dfrac{f(t)-f(1)}{t-1}, & t \neq 1 \\ f'(1), & t=1, \end{cases}$ where $f:[0,1] \to \Re$ is a Riemann integrable function having a finite derivative at $t = 1$, and $\omega(t) := 1-t$. By use of the identity $\int_0^1 f(x)dx \equiv f(1) - \int_0^1 f_\omega(x)\omega(x)dx$ and slight algebraic manipulation in the standard formulas [4, (1.4.7), (3.1.13), (3.1.26)], we express the $n$-point Gauss-type rules as follows:

(a) Gauss Rule:
$$\int_0^1 f(x)dx \approx f(1) - \sum_{i=1}^n v_i^G f_\omega(t_i^G), \tag{1}$$

(b) Gauss-Radau Rule (Right end):
$$\int_0^1 f(x)dx \approx f(1) - \sum_{i=1}^n v_i^G f_\omega(t_i^G), \tag{2}$$

(c) Gauss-Lobatto Rule
$$\int_0^1 f(x)dx \approx f(1) + v_1^L(f(0)-f(1)) - \sum_{i=2}^{n-2} v_i^L f_\omega(t_i^L) - v_n^L f'(1). \tag{3}$$

Here, $(t_i^G, v_i^G), (t_i^R, v_i^R)$ and $(t_i^L, v_i^L)$ respectively are the pairs of nodes and weights for the respective rules (a), (b) and (c). These nodes and weights arise from the zeros of polynomials [3,4] which are orthogonal with respect to weight function $\omega(t)$ over $[0, 1]$. If we denote such polynomials by $p_k$ with $p_{-1}(t)=0, p_0(t)=1$, then $p_k, k=0,1,2,\ldots$ are generated from the three-term recurrence relation [4]:
$$p_{k+1}(t) = (t-\alpha_{k+1})p_k(t) - \beta_{k+1}p_{k-1}(t) \tag{4}$$
where the recursion coefficients are given by
$$\alpha_{k+1} = \frac{\langle tp_k, p_k \rangle_\omega}{\langle p_k, p_k \rangle_\omega}, \quad \beta_{k+2} = \frac{\langle p_{k+1}, p_{k+1} \rangle_\omega}{\langle p_k, p_k \rangle_\omega}; \quad k=0,1,2,\ldots. \tag{5}$$

with $\beta_1 = 0$. Here, the notation $\langle F, G \rangle_\omega$ stands for $\int_0^1 F(t)G(t)\omega(t)dt$. The knowledge of $\alpha_k$ and $\beta_k$ in (5) allows the zeros of $p_k, k=0,1,2,\ldots$ and the weights to be readily calculated from the eigenvalues and eigenvectors of a symmetric tri-diagonal matrix [7, 9, 13]. In our context, a formal procedure as suggested by Golub and Welsh (see [7, 9]) is summarized in

**Theorem A.** *With $\alpha_n$ and $\beta_n$ as defined in (5), the "n free nodes" of the Gauss-type quadrature formulas with respect to weight function $\omega(t)$ over the interval $[0,1]$ are precisely the eigenvalues of the Jacobi matrix*

$$J_{n,\omega} = \begin{bmatrix} \alpha_1 & \sqrt{\beta_2} & 0 & & \cdots & & 0 \\ \sqrt{\beta_2} & \alpha_2 & \sqrt{\beta_3} & 0 & & & 0 \\ 0 & \sqrt{\beta_3} & \alpha_3 & & & & \\ & & & \ddots & \ddots & \ddots & \vdots \\ & & & \ddots & \alpha_n & \sqrt{\beta_{n+1}} & 0 \\ \vdots & & & \ddots & \sqrt{\beta_{n+1}} & \alpha^R_{n+1} & \sqrt{\beta^L_{n+2}} \\ 0 & 0 & & \cdots & 0 & \sqrt{\beta^L_{n+2}} & \alpha^L_{n+2} \end{bmatrix}, \qquad (6)$$

*where in case of*

(a) *Gauss rule:* $J_{n,\omega}$ *is considered after deleting the last two columns and the last two rows.*

(b) *Gauss-Radau rule (Right end): the last column and the last row in $J_{n,\omega}$ are deleted. In addition we set $\alpha^R_{n+1} := 1 - \beta_{n+1}\dfrac{p_{n-1}(1)}{p_n(1)}$.*

(c) *Gauss-Lobatto rule: we choose $\alpha^R_{n+1} = \alpha_{n+1}$ in $J_{n,\omega}$ whereas the entries $\alpha^L_{n+2}$ and $\beta^L_{n+2}$ are the solution of the following 2x2 linear system:*

$$\begin{bmatrix} p_{n+1}(0) & p_n(0) \\ p_{n+1}(1) & p_n(1) \end{bmatrix} \begin{bmatrix} \alpha^L_{n+2} \\ \beta^L_{n+2} \end{bmatrix} = \begin{bmatrix} 0 \\ p_{n+1}(1) \end{bmatrix}.$$

Moreover, the respective weights in (1)-(3) are given by $wu^2_{i,1}, i = 1,2,3,\ldots$, where $w = \int_0^1 \omega(x)dx$ and $u_{i,1}$ is the first component of the associated normalized eigenvector $u_i$.

**Remark 1:** The quadrature rules (a)-(c) are applicable to appropriate integrands $g(x)$ over a finite interval $[a,b]$ by considering $\int_a^b g(x)dx = \int_0^1 f(t)dt$ with $f(t) := (b-a)g((b-a)t + a), t \in [0,1]$. In addition, the fixed node "1" in the rules (1)-(3), if required, can be replaced by 0, the left end point of the interval, by expressing $\int_0^1 g(x)dx = \int_0^1 f(t)dt$ where $f(t) := g(1-t), t \in [0,1]$.

We are also interested to observe accumulation of the round-off error that may result while programming the proposed rules. Based on the structure of matrix (6), the computational efficacy, indeed, relies on appropriate representation of recursion coefficients $\alpha_n, \beta_n$ (cf. (5)). The identity-type polynomials discussed below provide one such representation.

## Identity-type Polynomials

It is known that the identity-type function [1, (2.16)]

$$\widehat{e}(t;c) := \frac{\Gamma(c)}{\Gamma(c-1)} \sum_{m=0}^{\infty} \frac{(c)_m (-c)_m}{(m!)^2} t^m, \tag{7}$$

with $(c)_0 = 1$ and $(c)_n = c(c-1)(c-2)....(c-n+1)$ satisfies the second order differential equation $t(1-t)\frac{d^2 y}{dt^2} + (1-t)\frac{dy}{dt} + c^2 y = 0, \ c > 0.$

Since $(n)_m = 0$ for $m \geq n+1$, the hypergeometric series $\sum_{m=0}^{\infty} \frac{(c)_m (-c)_m}{(m!)^2} t^m$ in (7) with c = n results to an nth degree polynomial

$$e_n(t) := \sum_{m=0}^{\infty} \frac{(n)_m (-n)_m}{(m!)^2} t^m, \ n = 1, 2, 3, ......, \tag{8}$$

which is the solution of the Sturm-Liouville problem [1, (4.1)]

$\frac{d}{dt}\left(t\frac{dy}{dt}\right) + \frac{n^2}{1-t} y = 0$. To meet our objective, we reproduce some properties of $e_n$ as discussed in [1]:

(i) The polynomials $e_n, n = 1, 2, ...$, are orthogonal [1, (4.3)] with respect to the weight function $\omega^*(t) = \frac{1}{(1-t)}$ over $[0,1]$. In fact,

$$\langle e_n, e_m \rangle_{\omega^*} = \begin{cases} 0 & \text{if } n \neq m, \\ \frac{1}{2n} & \text{if } n = m. \end{cases} \tag{9}$$

(ii) $(1-t)$ is a factor of each $e_n$, n = 1, 2,... [1, (3.1)]. In particular, if we write

$$e_n(t) := (1-t) e_{n-1}^*(t), \tag{10}$$

then the first few factor polynomials $e_{n-1}^*$, n = 1,2,..., are given by:

$e_0^*(t) = 1,$

$e_1^*(t) = (1-3t),$

$e_2^*(t) = (1-8t+10t^2),$

$\vdots$

The polynomials $e_{n-1}^*$, n = 1,2,3,... ., are orthogonal with respect to $\omega(t) := 1-t$ over $[0,1]$ since (cf (9)-(10))

$$\langle e_{n-1}^*, e_{m-1}^* \rangle_{\omega} = \langle e_n, e_m \rangle_{\omega^*}, \ n, m = 1,2,3,.... \tag{11}$$

Thus, letting $\kappa_n$ be the leading coefficient of $e_n$, it follows that

$$p_{n-1}(t) := \frac{-1}{\kappa_n} e_{n-1}^*(t), \quad n = 1, 2,...... \tag{12}$$

are monic as well as orthogonal polynomials over $[0,1]$ with respect to $\omega(t)$. Here, we observe that $\omega(t)$ as defined above is the weight function required in the proposed quadrature formulas (1)-(3).

## Representation of Recursion Coefficients for $p_n$

In view of (6), we are interested in simple representation for the recursion coefficients (cf (5)) that generate the orthogonal polynomials $p_{n-1}$ (cf (12)) via relation (5). These representations are obtained by invoking some properties of the factor polynomials $e_{n-1}^*$.

**Lemma 1.** Let $p_{n-1}$ be the orthogonal polynomial and $\kappa_n$ the constant as described in (12). Then

$$\kappa_n = (-1)^n \frac{(2n-1)!}{n!(n-1)!}, \qquad (13)$$

$$p_{n-1}(0) = -\kappa_n^{-1}, \qquad (14)$$

$$p_{n-1}(1) = (-1)^n n \kappa_n^{-1}. \qquad (15)$$

**Proof.** Note that $\kappa_n$, the leading coefficient of $e_n$ can be expressed (cf (8)) as $\kappa_n = \frac{(n)_n (-n)_n}{(n!)^2}$ where $(-n)_n = (-1)^n n(n+1)\ldots\ldots(2n-1)$. This leads us to (13).

Next, observe that (cf (10) and (8)) $e_{n-1}^*(0) = e_n(0) = 1$. Thus, from (12), after replacing $t$ by 0 we obtain (14).

The proof of relation (15) requires the notion of Gauss hypergeometric functions and some transformations as discussed in [1]. Using [1, (3.5)], it follows that $e_{n-1}^*(1) = (-1)^{n+1} n$. Substituting this in (10) leads to (15) which completes the proof of Lemma1.

Our main result is as follows:

**Theorem 2.** *An explicit representation of $\alpha_n$ and $\beta_{n+1}$ in Jacobi matrix (6) is given by:*

$$\beta_{n+1} = \frac{n(n+1)}{4(2n+1)^2}, \qquad (16)$$

$$\alpha_n = \frac{2n^2 - 1}{4n^2 - 1}. \qquad (17)$$

*Moreover, the additional entries in (6) which are required respectively for Gauss-Radau and Gauss-Lobatto rules (cf (2)-(3)) are expressed as*

$$\alpha_{n+1}^{R,1} = \frac{3n^2 + 6n + 2}{4n^2 + 6n + 2}, \qquad (18)$$

*and*

$$\left.\begin{array}{l}\alpha_{n+2}^{L} = \dfrac{n+2}{2n+3} \\[6pt] \beta_{n+2}^{L} = \dfrac{(n+2)^{2}}{2(2n+3)^{2}}\end{array}\right\}. \qquad (19)$$

**Proof.** Recall that $\alpha_n$ and $\beta_n$ in (6) are the recursion coefficients defined in (5) which generate the orthogonal polynomials $p_{n-1}(t)$ as indicated in (4). To prove (16), we use (12) in the expression given for $\beta_n$ in (5) and then apply (11) to get

$\beta_{n+1} = \left(\dfrac{\kappa_n}{\kappa_{n+1}}\right)^2 \dfrac{\langle e_{n+1}, e_{n+1}\rangle_{\omega^*}}{\langle e_n, e_n\rangle_{\omega^*}} = \left(\dfrac{\kappa_n}{\kappa_{n+1}}\right)^2 \dfrac{n}{n+1}$. From (13), note that

$$\dfrac{\kappa_n}{\kappa_{n+1}} = \dfrac{-(n+1)}{2(2n+1)}. \qquad (20)$$

Inserting this in the preceding expression provides (16).

To justify (17), set $t = 0$ and $k = n - 1$ in equation (4). This gives us

$$-1 = \alpha_n \dfrac{\kappa_{n+1}}{\kappa_n} + \beta_n \dfrac{\kappa_{n+1}}{\kappa_{n-1}}. \qquad (21)$$

The second quotient on the right side of (21) can be expressed as $\dfrac{\kappa_{n+1}}{\kappa_{n-1}} = \dfrac{4(n^2-1)}{n(n+1)}$ (cf (13)). Thus, by (16), $\beta_n \dfrac{\kappa_{n+1}}{\kappa_{n-1}} = \dfrac{(2n+1)(n-1)}{(2n-1)(n+1)}$. Using (20) and this expression in (21) we arrive at (17).

To prove (18), we consider the relation given in Theorem A(b): $\alpha_{n+1}^{R} = 1 - \beta_{n+1}\dfrac{p_{n-1}(1)}{p_n(1)}$. After appropriate substitutions from (13)-(16), we note that $\beta_{n+1}\dfrac{p_{n-1}(1)}{p_n(1)}$ reduces to $\dfrac{n^2}{2(n+1)(2n+1)}$. Thus, $\alpha_{n+1}^{R} = \dfrac{3n^2+6n+2}{4n^2+6n+2}$ which is the desired relation (18).

For (19), we revisit Theorem A(c) and note from the solution of 2X2 system that

$$\alpha_{n+2}^{L} = \dfrac{-p_{n+1}(1)p_n(0)}{D}, \quad \beta_{n+2}^{L} = \dfrac{-p_{n+1}(0)p_{n+1}(1)}{D} \qquad (22)$$

with $D = p_{n+1}(0)p_n(1) - p_{n+1}(1)p_n(0)$. Now applying the relations (13)-(15) to (22) and then performing some algebraic manipulations, we obtain (19). This completes the proof.

**Remark 2:** Theorem 2 provides simple expressions for $\alpha_n, \beta_n$. Thus the nodes and weights required in the proposed quadrature formulae (1) - (3) are easily obtainable by an application of Theorem A.

## Numerical Examples

We have applied the quadrature rules (1)-(3) to four functions having different characteristics:

(1) $f_1(x) = x^{1/8}$, (2) $f_2(x) = \cosh^2(5(x-0.5))/5$, (3) $f_4(x) = \dfrac{2}{2+\sin(10\pi x)}$,

(4) $f_3(x) = \dfrac{1}{(x-0.3)^2 + 0.01} + \dfrac{1}{(x-0.9)^2 + 0.04}$.

Similar functions have been considered as test integrands in [2, 11, 12]. We have computed percentage errors that occurred in the approximation of all four integrals by using *n*-Gauss-type rules, *n* = 2, 3,…, 11, with respect to both weight functions $\omega(t) := 1-t$ and $\omega(t) := 1$. The data for each integrand is provided in separate tables for the sake of comparison. The following abbreviations have been used to identify the rule in each table:

- G(1) := Gauss rule w.r.t. $\omega(t) := 1$,
- G (1 – *t*) := Gauss rule w.r.t. $\omega(t) := 1-t$,
- G-R(1) := Gauss-Radau rule (Right end point) w.r.t. $\omega(t) := 1$,
- G-R(1 – *t*) := Gauss-Radau rule (Right end point) w.r.t. $\omega(t) := 1-t$,
- G-L(1) := Gauss-Lobatto rule w.r.t. $\omega(t) := 1$,
- G-L(1 – *t*) := Gauss-Lobatto rule w.r.t. $\omega(t) := 1-t$.

**Fig 1: Graph of function $f_1$**

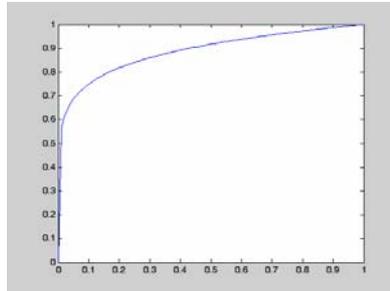

**Table 1**: Comparison of Percentage Errors for $f_1$

| n | %Error G (1) | %Error G (1-*t*) | %Error G-R (1) | %Error G-R (1-*t*) | %Error G-L (1) | %Error G-L (1-*t*) |
|---|---|---|---|---|---|---|
| 2 | 0. 922 | **0.622** | 0.622 | **0.120** | 5.693 | **4.424** |
| 3 | 0. 423 | **0.317** | 0.317 | **0.041** | 3.198 | **2.604** |
| 4 | 0. 238 | **0.189** | 0.189 | **0.063** | 2.025 | **1.702** |
| 5 | 0. 151 | **0.125** | 0.125 | **0.058** | 1.386 | **1.192** |
| 6 | 0.103 | **0.088** | 0.088 | **0.049** | 1.002 | **0.878** |
| 7 | 0. 075 | **0.065** | 0.065 | **0.041** | 0.755 | **0.671** |
| 8 | 0.056 | **0.049** | 0.049 | **0.034** | 0.587 | **0.527** |
| 9 | 0.043 | **0.039** | 0.039 | **0.028** | 0.468 | **0.425** |
| 10 | 0.035 | **0.031** | 0.031 | **0.024** | 0.381 | **0.349** |
| 11 | 0.028 | **0.026** | 0.026 | **0.020** | 0.316 | **0.291** |

### Fig 2: Graph of function $f_2$

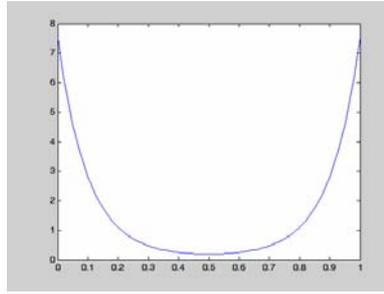

### Table 2: Comparison of Percentage Errors for $f_2$

| $n$ | %Error G (1) | %Error G (1-$t$) | | %Error G-R (1) | %Error G-R (1-$t$) | | %Error G-L (1) | %Error G-L (1-$t$) |
|---|---|---|---|---|---|---|---|---|
| 2 | 36.901 | **3.082** | | 3.082 | **14.037** | | 9.284 | **1.078** |
| 3 | 6.523 | **0.202** | | 0.202 | **1.103** | | 0.834 | **0.045** |
| 4 | 0.647 | **$9.72 \times 10^{-3}$** | | $9.72 \times 10^{-3}$ | **0.06.2** | | 0.050 | **$1.49 \times 10^{-3}$** |
| 5 | 0.041 | **$3.48 \times 10^{-4}$** | | $3.48 \times 10^{-4}$ | **$2.57 \times 10^{-3}$** | | $2.16 \times 10^{-3}$ | **$3.87 \times 10^{-5}$** |
| 6 | $1.83 \times 10^{-3}$ | **$9.54 \times 10^{-6}$** | | $9.54 \times 10^{-6}$ | **$8.00 \times 10^{-5}$** | | $6.92 \times 10^{-5}$ | **$8.02 \times 10^{-7}$** |
| 7 | $6.02 \times 10^{-5}$ | **$2.05 \times 10^{-7}$** | | $2.05 \times 10^{-7}$ | **$1.93 \times 10^{-6}$** | | $1.70 \times 10^{-6}$ | **$1.35 \times 10^{-8}$** |
| 8 | $1.50 \times 10^{-6}$ | **$3.58 \times 10^{-9}$** | | $3.58 \times 10^{-9}$ | **$3.72 \times 10^{-8}$** | | $3.33 \times 10^{-8}$ | **$1.88 \times 10^{-10}$** |
| 9 | $2.99 \times 10^{-8}$ | **$5.15 \times 10^{-11}$** | | $5.12 \times 10^{-11}$ | **$5.84 \times 10^{-10}$** | | $5.28 \times 10^{-10}$ | **$2.13 \times 10^{-12}$** |
| 10 | $4.79 \times 10^{-10}$ | **$4.76 \times 10^{-13}$** | | $6.86 \times 10^{-13}$ | **$7.80 \times 10^{-12}$** | | $6.95 \times 10^{-12}$ | **$2.80 \times 10^{-14}$** |
| 11 | $6.29 \times 10^{-12}$ | **$1.40 \times 10^{-14}$** | | $8.41 \times 10^{-14}$ | **$3.36 \times 10^{-13}$** | | $9.81 \times 10^{-14}$ | **$1.26 \times 10^{-13}$** |

### Fig 3: Graph of function $f_3$

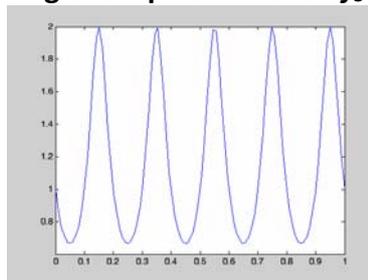

**Table 3**: Comparison of Percentage Errors for $f_3$

| n | %Error G (1) | %Error G (1-$t$) | %Error G-R (1) | %Error G-R (1-$t$) | %Error G-L (1) | %Error G-L (1-$t$) |
|---|---|---|---|---|---|---|
| 2 | 10.688 | **3.721** | 3.721 | **58.697** | 4.105 | **20.742** |
| 3 | 11.511 | **13.906** | 13.906 | **21.056** | 5.516 | **27.428** |
| 4 | 3.733 | **20.939** | 20.939 | **10.518** | 3.904 | **0.387** |
| 5 | 1.804 | **23.850** | 23.850 | **12.648** | 3.032 | **35.580** |
| 6 | 0.0102 | **49.801** | 49.801 | **1.632** | 8.408 | **36.267** |
| 7 | 7.133 | **22.412** | 22.412 | **10.528** | 0.486 | **11.073** |
| 8 | 4.733 | **12.604** | 12.604 | **1.417** | 3.160 | **4.061** |
| 9 | 4.080 | **3.055** | 3.055 | **1.269** | 1.214 | **0.605** |
| 10 | 0.585 | **1.059** | 1.059 | **4.751** | 5.516 | **0.150** |
| 11 | 5.71 | **0.535** | 0.535 | **3.091** | 2.617 | **1.017** |

**Fig 4: Graph of function $f_4$**

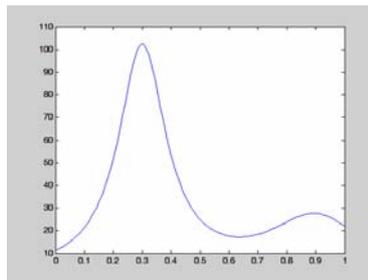

**Table 4**: Comparison of Percentage Errors for $f_4$

| n | %Error G (1) | %Error G (1-$t$) | %Error G-R (1) | %Error G-R (1-$t$) | %Error G-L (1) | %Error G-L (1-$t$) |
|---|---|---|---|---|---|---|
| 2 | 12.988 | **32.88** | 32.884 | **34.468** | 42.957 | **0.143** |
| 3 | 29.185 | **9.812** | 9.812 | **17.083** | 20.888 | **17.980** |
| 4 | 23.274 | **17.806** | 17.806 | **2.339** | 4.551 | **15.027** |
| 5 | 5.735 | **12.173** | 12.173 | **7.241** | 4.340 | **4.979** |
| 6 | 5.474 | **4.624** | 4.624 | **7.420** | 6.547 | **1.502** |
| 7 | 6.900 | **0.919** | 0.919 | **3.496** | 4.076 | **3.257** |
| 8 | 3.740 | **2.995** | 2.995 | **0.175** | 0.792 | **2.344** |
| 9 | 0.567 | **2.323** | 2.323 | **1.265** | 1.018 | **0.792** |
| 10 | 1.038 | **0.806** | 0.806 | **1.244** | 1.229 | **0.268** |
| 11 | 1.195 | **0.249** | 0.249 | **0.609** | 0.672 | **0.581** |

## Conclusion

We have introduced Gauss-type quadrature rules (cf. (1) - (3)) on the interval [0, 1] which, except for rule (1), involve derivative of the integrand at the right end of the interval. These rules also preserve the exactness and convergence properties like the corresponding Gauss quadrature rules. Our proposed rules also utilize an additional functional value as compared to the rules based on weight function 1. Four functions of different characteristic were selected for implementation of the proposed rules. Percentage error is used as a criterion of comparative accuracy.

For $f_1$ we note that the rules using weight function $\omega(t) := 1-t$ are consistently better than those for weight function $\omega(t) := 1$. Further, Gauss-Radau is the best rule for this function and Gauss-Lobatto the worst. The difference between the accuracy of the Gauss scheme and the Gauss-Radau scheme is relatively small. Since the two schemes for each of the rules (1)-(3) are of nearly the same percentage accuracy, it is worthwhile to introduce a measure of *relative accuracy*, $\alpha_r$, defined as the ratio of the percentage accuracies. From Table 1.1 related to the Gauss scheme, we note that $\alpha_r \approx 1.13$ at $n = 10$ in favor of the weight function $\omega(t) := 1-t$. At $n = 35$, it reduces to about 1.03. For the Gauss-Radau scheme at $n = 10$, we have $\alpha_r \approx 1.29$. At $n = 35$ it reduces to about 1.05. For the Gauss-Lobatto rule, the relative accuracy is about the same as observed for the Gauss scheme. Thus, for smaller $n$ the differences are significant. It is interesting to note that the percentage accuracy of the Gauss scheme with $\omega(t) := 1-t$ is about the same as that of the Gauss-Radau with weight factor $\omega(t) := 1$ for all $n$ that we checked for the function $f_1$.

For $f_2$ the percentage error is extremely small at $n = 10$ in all cases. Checking at $n = 20$ we see that it has already reached the limit of numerical stability at $n = 10$ or so. Here the Gauss-Lobatto is the best and the Gauss-Radau the worst. Further, the weight function $\omega(t) := 1$ is *better* for the Gauss-Radau than $\omega(t) := 1-t$. Here $\alpha_r \approx 11.4$ in favour of $\omega(t) := 1$. On the other hand for the Gauss scheme, $\alpha_r \approx 993.7$ in favour of $\omega(t) := 1-t$. This is the biggest advantage for this weight function. For the Gauss-Lobatto we have $\alpha_r \approx 40.3$ in favour of $\omega(t) := 1-t$.

The test function, $f_3$, is highly oscillating. Here, we note that the percentage error fluctuates wildly with changes in $n$ for all the six cases (cf. Tables 4.1 - 4.3). For the hump function ranges of 5 are seen to be adequate (cf Table 5). On the other hand, we noted that for the function $f_4$ ranges of 10 turn out to be barely adequate for larger values of $n$.

In the case of $f_4$, the famous humps function, the percentage error is seen to fluctuate up to $n = 11$. As such, we cannot draw any conclusions from it. One can look at the general change in accuracy as we determine average over each of the ranges: $n$ =12–16; 17 – 21; 22 – 26; 27 – 31; 32 – 36. We see that the error decreases as one goes to higher range. The differences are given in Table 5.

**Table 5: Average Accuracy of the Rules for Test Function $f_4$ over Various Ranges**

| Range | G(1) | G(1-t) | G-R(1) | G-R(1-t) | G-L(1) | G-L(1-t) |
|---|---|---|---|---|---|---|
| (12 – 16) | $2.5\times10^{-1}$ | $2.5\times10^{-1}$ | $2.5\times10^{-1}$ | $1.2\times10^{-1}$ | $1.3\times10^{-1}$ | $1.6\times10^{-1}$ |
| (17 – 21) | $2.1\times10^{-2}$ | $2.7\times10^{-2}$ | $2.7\times10^{-2}$ | $2.2\times10^{-2}$ | $2.1\times10^{-2}$ | $1.3\times10^{-2}$ |
| (22 -26) | $3.9\times10^{-3}$ | $2.3\times10^{-3}$ | $2.3\times10^{-3}$ | $2.6\times10^{-3}$ | $2.6\times10^{-3}$ | $1.8\times10^{-3}$ |
| (27 –31) | $4.6\times10^{-4}$ | $3.3\times10^{-4}$ | $3.3\times10^{-4}$ | $2.3\times10^{-4}$ | $2.3\times10^{-4}$ | $2.6\times10^{-4}$ |
| (32-26) | $4.0\times10^{-5}$ | $4.4\times10^{-5}$ | $4.4\times10^{-5}$ | $2.6\times10^{-5}$ | $2.6\times10^{-5}$ | $2.4\times10^{-5}$ |

We note that the identity between G(1-t) and G-R(1) carries over in this case despite of the fluctuations. Further, for this function, G(1), G(1–t), G-R(1) are comparable in accuracy and the same observation holds for G-R(1-t), G-L(1) and G-L(1–t). The latter group is more accurate than the former group.

The test function $f_3$ is highly oscillating. Here, we note that the percentage error fluctuates wildly with changes in $n$ for all the six cases (cf Tables 3). For the hump function ranges of 5 are seen to be adequate (cf Table 5). On the other hand, we noted that for the function $f_3$ ranges of 10 turn out to be barely adequate for larger values of $n$.

It is worthwhile to point out the stability of the proposed methods. We have observed stability of the three methods up to $n$ = 2000 in case of $\omega(t):=1-t$. On the other hand, the methods G-L(1) and G-R(1) fail to work respectively after $n$ = 272 and $n$ = 539 because of the following observations in the MATLAB programming:

i) The 2X2 linear system in Theorem A(c) which determines the additional parameters for G-L(1) method turns out inconsistent beyond $n$ = 272.

ii) In case of G-R(1), the value $p_n(1)$ after $n$ = 539 becomes zero.

As mentioned in Remark 1 the proposed rules can be easily used for any finite interval $[a,b]$ instead of $[0,1]$ by appropriate translation and re-scaling.

## Acknowledgment

Both authors are grateful to KFUPM for the excellent research facilities availed during the preparation of this article.